\documentclass{article}

\title{Contiguous relations of hypergeometric series}

\author{Raimundas Vid\=unas\footnote{Supported by NWO,
project number 613-06-565.}\\ \em University of Amsterdam}

\newtheorem{theorem}{Theorem}[section]

\newtheorem{lemma}[theorem]{Lemma}
\newcommand{\proof}{{\bf Proof. }}
\newcommand{\QED}{\hfill \mbox{QED.}\\}

\newcommand{\hpg}[5]{{}_{#1}\mbox{\rm F}_{\!#2}\!
  \left(\left.{#3 \atop #4}\right|\, #5 \right) }

\newcommand{\ghpg}[2]{{{\bf F}\!\left(#1 \atop #2\right)}}
\newcommand{\hpgo}[2]{{}_{#1}\mbox{\rm F}_{\!#2}}

\newcommand{\qhpg}[5]{{}_{#1}\phi_{#2}\!
  \left(\left.{#3 \atop #4}\,\right|q;\,#5 \right) }

\newcommand{\qhpgo}[2]{{}_{#1}\phi_{#2}}

\newcommand{\shif}[1]{{\cal S}_#1}
\newcommand{\bfo}[1]{{\bf #1}}
\newcommand{\emphi}[1]{\widetilde{#1}}

\newcommand{\CC}{\mbox{\bf C}}

\begin{document}

\date{}
\maketitle

\begin{abstract}
The 15 Gauss contiguous relations for ${}_2\mbox{F}_1$ hypergeometric series
imply that any three ${}_2\mbox{F}_1$ series whose corresponding parameters
differ by integers are linearly related (over the field of rational
functions in the parameters). We prove several properties of coefficients of
these general contiguous relations, and use the results to propose effective
ways to compute contiguous relations. We also discuss contiguous relations
of generalized and basic hypergeometric functions, and several
applications of them.\\
{\em Keywords:} Hypergeometric function, contiguous relations.
\end{abstract}

\section{Contiguous relations of $\hpgo{2}{1}$ series}

In this paper, let
\begin{equation}
\ghpg{a,\,b}{c}:=\hpg{2}{1}{a,\,b}{c}{z}=1+\frac{a\,b}{c\cdot 1}\,z+
\frac{a\,(a+1)\,b\,(b+1)}{c\,(c+1)\cdot 1\cdot 2}\,z^2+\ldots
\end{equation}
denote the Gauss hypergeometric function with the argument $z$.
Two hypergeometric functions with the same argument $z$ are {\em contiguous}
if their parameters $a$, $b$ and $c$ differ by integers. For example,
$\ghpg{a,\,b}{c}$ and $\ghpg{a+10,\,b-7}{c+3}$ are contiguous. As is known
\cite[2.5]{specfaar}, for any three contiguous $\hpgo{2}{1}$ functions there
is a {\em contiguous} relation, which is a linear relation, with coefficients
being rational functions in the parameters $a$, $b$, $c$ and the argument $z$.
For example,
\begin{eqnarray} \label{contigaa}
a\,(z-1)\,\ghpg{a+1,\,b}{c}+(2a-c-az+bz)\,\ghpg{a,\,b}{c} \nonumber\\
+(c-a)\,\ghpg{a-1,\,b}{c} & = & 0, \\
%\label{contigbb} c\,(1-z)\,\ghpg{a,\,b}{c}-c\,\ghpg{a-1,\,b}{c}
%+(c-b)\,z\;\ghpg{a,\,b}{c+1} & = & 0,\\
\label{contigac} a\,\ghpg{a+1,\,b}{c}-(c-1)\,\ghpg{a,\,b}{c-1}
+(c-a-1)\,\ghpg{a,\,b}{c} & = & 0,\\
\label{contigab} a\,\ghpg{a+1,\,b}{c}-b\,\ghpg{a,\,b+1}{c}
+(b-a)\,\ghpg{a,\,b}{c} & = & 0.
\end{eqnarray}
In these relations two hypergeometric series differ just in one
parameter from the third hypergeometric series, and the difference
is 1. The relations of this kind were found by Gauss; there are 15
of them. A contiguous relation between any three contiguous
hypergeometric functions can be found by combining linearly a
sequence of Gauss contiguous relations. In the next section we
discuss this and other ways to compute contiguous relations.

The following theorem summarizes some properties of coefficients
of contiguous relations. These results are useful in computations
and applications of contiguous relations. We assume that the
parameters $a$, $b$, $c$ and $z$ are not related, and by
$\shif{a}$, $\shif{b}$, $\shif{c}$ we denote the shift operators
$a\mapsto a+1$, $b\mapsto b+1$ and $c\mapsto c+1$ respectively.
%The theorems of this section are proved in section \ref{thproofs}.
\begin{theorem} \label{spcontig}
For any integers $k$, $l$, $m$ there are unique functions
$\bfo{P}(k,l,m)$ and $\bfo{Q}(k,l,m)$, rational in the parameters
$a$,$b$,$c$ and $z$, such that
\begin{equation} \label{contklm}
\ghpg{a+k,\,b+l}{c+m} = \bfo{P}(k,l,m)\;\ghpg{a,\,b}{c}
+\bfo{Q}(k,l,m)\,\ghpg{a+1,\,b}{c}.
\end{equation}
These functions satisfy the same contiguous relations as
$\ghpg{a,\,b}{c}$, that is, if %{\rm (\ref{gencontig})}
\begin{equation} \label{gencontig}
\sum_{n=1}^3 f_n\;\ghpg{a+\alpha_n,\,b+\beta_n}{c+\gamma_n}=0,
\qquad f_1,f_2,f_3\in\CC(a,b,c,z),
\end{equation}
is a contiguous relation, then $\sum_{n=1}^3 \left(
\shif{a}^k\,\shif{b}^l\,\shif{c}^m\,f_n \right)\;
\bfo{P}(k\!+\!\alpha_n,\,l\!+\!\beta_n,\,m\!+\!\gamma_n)=0$, and
similarly for $\bfo{Q}(k,l,m)$.
%$\sum_{n=1}^3 (\shif{a}^k\,\shif{b}^l\,\shif{c}^m\,f_n)\;
%\bfo{Q}(k\!+\!\alpha_n,\,l\!+\!\beta_n,\,m\!+\!\gamma_n)=0$.

Besides, the following expressions hold:
\begin{equation} \label{initcond}
\begin{array}{ccc}
\bfo{P}(0,0,0)=1, & & \bfo{Q}(0,0,0)=0,\\
\bfo{P}(1,0,0)=0, & & \bfo{Q}(1,0,0)=1,
\end{array}
\end{equation}
\begin{equation} \label{pinq} \bfo{P}(k,\,l,\,m) =
\frac{c-a-1}{(a+1)\,(1-z)}\;\shif{a}\bfo{Q}(k-1,\,l,\,m),
\end{equation}
\begin{eqnarray} \label{qargplus}
\bfo{Q}(k+k',\,l+l',\,m+m')\ =\ \left( \shif{a}^k \shif{b}^l
\shif{c}^m\, \bfo{Q}(k',l',m') \right)\,\bfo{Q}(k\!+\!1,l,m) \nonumber\\
+\left( \shif{a}^{k+1}\!\shif{b}^l \shif{c}^m \frac{c-a}{a\,(1-z)}
\bfo{Q}(k'\!-\!1,l',m') \right) \bfo{Q}(k,l,m), \hspace{-4pt}
\end{eqnarray}
\begin{eqnarray} \label{pqdet} \left|\begin{array}{cc}
\bfo{P}(k,l,m) & \bfo{P}(k',l',m') \\
\bfo{Q}(k,l,m) & \bfo{Q}(k',l',m') \end{array}\right| & \!\!=\!\!
&  \frac{(c)_m\;(c)_m\;z^{-m}\,(z-1)^{m-k-l}}
{(a\!+\!1)_k\,(b)_l\,(c\!-\!a)_{m-k}(c\!-\!b)_{m-l}} \nonumber \\
& & \hspace{-2pt} \times \shif{a}^k\,\shif{b}^l\,\shif{c}^m\,
\bfo{Q}(k'-k,l'-l,m'-m),
\end{eqnarray}
\begin{eqnarray} \label{minklm} \bfo{Q}(-k,-l,-m) &\!\!=\!\!&
\frac{(-1)^{m+1}\;(-a)_k\;(1\!-\!b)_l\;z^m\,(1-z)^{k+l-m}}
{(1\!-\!c)_m\;(1\!-\!c)_m\;(c\!-\!a)_{k-m}\;(c\!-\!b)_{l-m}} \nonumber\\
& & \times
\shif{a}^{-k}\,\shif{b}^{-l}\,\shif{c}^{-m}\,\bfo{Q}(k,l,m).
\end{eqnarray}
Here $(a)_k=\Gamma(a+k)/\Gamma(a)$ is the Pochhammer symbol; for positive $k$
it is equal to $a\,(a+1)\,\ldots\,(a+n-1)$.
% well defined for negative $k$ as well.
\end{theorem}
This theorem is proved in Section \ref{thproofs}. In the following section,
we discuss computational aspects and applications of contiguous relations,
including contiguous relations of generalized and basic hypergeometric
series.

\section{Computational aspects and applications}

\bfo{Computational aspects.} To compute a contiguous relation
(\ref{gencontig}) between three $\hpgo{2}{1}$ series, one can take
one of the series as $\ghpg{a,\,b}{c}$. If necessary, one can
apply a suitable shift operator to the contiguous relation
afterwards. Then the other two hypergeometric series are
$\ghpg{a+k',\,b+l'}{c+m'}$ and $\ghpg{a+k'',\,b+l''}{c+m''}$ for
some integers $k',l',m'$ and $k'',l'',m''$. According to Theorem
\ref{spcontig}, these two functions can be expressed linearly in
$\ghpg{a,\,b}{c}$ and $\ghpg{a+1,\,b}{c}$ as in (\ref{contklm}).
Elimination of $\ghpg{a+1,\,b}{c}$ gives:
\begin{eqnarray} \label{contig2klm}
\bfo{Q}(k'',l'',m'')\,\ghpg{a+k',\,b+l'}{c+m'}-\bfo{Q}(k',l',m')\,
\ghpg{a+k'',\,b+l''}{c+m''} =\nonumber \\
\big( \bfo{P}(k',l',m')\,\bfo{Q}(k'',l'',m'')-
\bfo{P}(k'',l'',m'') \bfo{Q}(k',l',m') \big)\,\ghpg{a,\,b}{c}.
\end{eqnarray}
Due to formula (\ref{pinq}) it is enough to be able to compute
$\bfo{Q}(k,l,m)$ for any integers $k$, $l$, $m$. For this one can
take a finite sequence of integer triples $(k_i,l_i,m_i)$ which
starts with $(0,0,0)$, $(1,0,0)$, ends with $(k,l,m)$, and any two
neighboring triples differ just in one component precisely by 1.
Then all $\bfo{Q}(k_i,l_i,m_i)$ including $\bfo{Q}(k,l,m)$ can be
computed consequently using the simplest Gauss contiguous
relations. A possible choice of the sequence is to take all
triples with intermediate integer values in the first component
between $(0,0,0)$ and $(k,0,0)$, then all ``intermediate triples"
between $(k,0,0)$ and $(k,l,0)$, and finally all triples between
$(k,l,0)$ and $(k,l,m)$.

Formula (\ref{qargplus}) allows us to compute $\bfo{Q}(k,l,m)$ by ``divide
and conquer" techniques, that is, by reducing the shift $(k,l,m)$
recursively by half (more or less) at each step; compare with
\cite[4.6.3]{knuthii}. A straightforward algorithm of this type is the
following: compute $\bfo{Q}(k,0,0)$ and $\bfo{Q}(k\!+\!1,0,0)$ using
(\ref{qargplus}) with shifts in the first parameter only (so there would be
$O(\log k)$ intermediate $\bf Q$-terms, perhaps linearly growing in size);
then use (\ref{contigab}) to compute $\bfo{Q}(k,1,0)$; and then alternate
similar application of (\ref{qargplus}) and Gauss contiguity relations to
compute $\bfo{Q}(k,l,0)$, $\bfo{Q}(k,l\!+\!1,0)$, $\bfo{Q}(k,l,1)$ and
$\bfo{Q}(k,l,m)$. Notice that most formulas of Theorem \ref{spcontig} cannot
be used when $a$, $b$ or $c$ are specialized, since the corresponding shift
operators would not be defined then.

Returning to the computation of a contiguous relation (\ref{gencontig}),
note that one can take any other of the given three functions in
(\ref{contig2klm}) as $\ghpg{a,\,b}{c}$. Proceeding in the same way one
would obtain a contiguous relation involving, say, $\bfo{Q}(-k',-l',-m')$
and $\bfo{Q}(k''-k',l''-l',m''-m')$. The two contiguous relations must be
the same up to a rational factor and corresponding shifts in the parameters.
Formulas (\ref{pqdet}) and (\ref{minklm}) are the explicit relations
(between the coefficients) implied by this fact. Notice that combination of
(\ref{contig2klm}) and (\ref{pqdet}) gives an expression of a contiguous
relation with coefficients linear in $\bfo{Q}$'s.

Contiguous relations (\ref{contklm}) and (\ref{contig2klm}) can be also
computed by combining Gauss contiguous relations themselves in similar ways:
after choosing a sequence of contiguous hypergeometric series ``connecting"
$\ghpg{a+k,\,b+l}{c+m}$ and $\ghpg{a,\,b}{c}$ so that neighboring series
differ just in one parameter by 1, or by using a similar formula to
(\ref{qargplus}). But this is double work compared to computing the
$\bfo{Q}$'s.

Takayama \cite{takayamcontig} showed how Gr\"obner basis can be used
to compute contiguous relalations. Paule \cite{paulecont} has shown that
contiguous relations can be computed by a generalized version of 
Zeilberger's algorithm. However, for large shifts $k$, $l$, $m$ this method
is not efficient, because the growing
degree (in the discrete parameters) of coefficients of contiguous relations
imply larger linear problems. On the other hand, recurrence relations for
hypergeometric functions can be computed as special cases of contiguous
relations by the methods described here, alternatively to Zeilberger's
algorithm. Similar remarks hold also for computation of contiguous relations
of generalized and basic hypergeometric series. On the website
\cite{koornvid} there is a link to a {\sf Maple} package for computing
contiguous relations
for $\hpgo{2}{1}$ series. % (Similar packages for $\hpgo{3}{2}(1)$ and
%$\qhpgo{2}{1}$ series are under construction.)
\\

\noindent \bfo{Generalization.} The generalized hypergeometric
series is defined as:
\begin{equation}
\hpg{r}{s}{a_1\,\ldots,a_r}{b_1,\ldots,b_s}{z}:=\sum_{k=0}^{\infty}
\frac{(a_1)_k\,\ldots\,(a_r)_k}{(b_1)_k\ldots (b_s)_k\,k!}{z^k}.
\end{equation}
Two $\hpgo{p}{q}$ series are {\em contiguous} if their corresponding upper
and lower parameters differ by integers. As for $\hpgo{2}{1}$ series, there
are linear relations between contiguous generalized hypergeometric series.
In general they relate $1+\max(p,q+1)$ hypergeometric series, see \cite[\S
48]{rainvill} and \cite[3.7]{specfaar}. In particular, relations
(\ref{contigac}) and (\ref{contigab}) also hold for hypergeometric series
with more upper and lower parameters. These relations allow to transform any
difference operator to a difference operator with different shifts in one
upper parameter only (by lowering all other upper parameters and raising the
lower ones). It follows that general contiguous relations for $\hpgo{p}{q}$
series are generated by the relations of type
(\ref{contigac})--(\ref{contigab}), and a recurrence relation with the
shifts in one upper parameter.

Like for $\hpgo{2}{1}$ series, contiguous relations for a class of
$\hpgo{p}{q}$ functions can be computed by linearly combining a sequence of
simplest contiguous relations. One can derive a formula analogous to
(\ref{contklm}) with at most $\max(p,q+1)$ fixed hypergeometric functions on
the right-hand side. The coefficients of the fixed hypergeometric functions
would satisfy the contiguity relations of $\hpgo{p}{q}$ functions, with
corresponding initial conditions like (\ref{initcond}). Those coefficients
usually are not all related by a formula like (\ref{pinq}), unless the
hypergeometric functions under consideration satisfy three-term contiguous
relations (say, $\hpgo{3}{2}(1)$ functions).

Similarly, there are contiguous relations for basic hypergeometric (or
$q$-hypergeometric) series; see \cite[10.9]{specfaar} for the definition of
$\qhpgo{r}{s}$ series. Two such series are contiguous if their corresponding
upper and lower parameters differ by a power of the base $q$. Moreover,
multiplicative $q$-shifts in the argument $z$ of these functions can also be
allowed, since the $q$-shift in the argument can be expressed in $q$-shifts
of the parameters, say:
\[
a\;\qhpg{2}{1}{a,\,b}{c}{qz}=(a-1)\;\qhpg{2}{1}{aq,\,b}
{b}{z}+\qhpg{2}{1}{a,\,b}{c}{z}.
%a\;\qhpg{r}{s}{a,\,\ldots}{\ldots}{qz}=(a-1)\;\qhpg{r}{s}{aq,\,\ldots}
%{\ldots}{z}+\qhpg{r}{s}{a,\,\ldots}{\ldots}{z}.
\]
For any three contiguous $\qhpgo{2}{1}$ series, where also $q$-shifts in the
argument $z$ are allowed, there is a contiguous relation. Allowing
$q$-shifts in the argument is natural, because many transformations of basic
hypergeometric series mix the parameters $a,b,c$ and the argument $z$, see
\cite[(10.10.1)]{specfaar}. A $q$-differential equation for $\qhpgo{r}{s}$
series can be interpreted as a contiguous relation in this more general
sense, since it can be seen as a $q$-difference equation.
\\

\noindent \bfo{Evaluation.} Contiguous relations can be used to
evaluate a hypergeometric function which is contiguous to a
hypergeometric series which can be satisfactorily evaluated. For
example, Kummer's identity \cite[Cor.~3.1.2]{specfaar}
\begin{equation} \label{kummer}
\hpg{2}{1}{a,\;b}{1+a-b}{\,-1} =
\frac{\Gamma(1+a-b)\;\Gamma(1+\frac{a}{2})}
{\Gamma(1+a)\;\Gamma(1+\frac{a}{2}-b)}.
\end{equation}
can be generalized to
\begin{equation} \label{gkummer}
\hpg{2}{1}{a+n,\,b}{a-b}{-1}= P(n)\,
\,\frac{\Gamma(a-b)\,\Gamma(\frac{a+1}{2})}
{\Gamma(a)\,\Gamma(\frac{a+1}{2}-b)}+Q(n)\,
\,\frac{\Gamma(a-b)\,\Gamma(\frac{a}{2})}
{\Gamma(a)\,\Gamma(\frac{a}{2}-b)}.
\end{equation}
Here $n$ is an integer, the factors $P(n)$ and $Q(n)$ can be
expressed for $n\ge 0$ as
\begin{eqnarray} \label{whippn}
P(n) & = & \frac{1}{2^{n+1}}\;
\hpg{3}{2}{-\frac{n}{2},\,-\frac{n+1}{2},\,\frac{a}{2}-b}
{\frac{1}{2},\;\frac{a}{2}}{\,1\,}, \\
\label{whipqn} Q(n) & = & \frac{n\!+\!1}{2^{n+1}}\;
\hpg{3}{2}{-\frac{n-1}{2},\,-\frac{n}{2},\,\frac{a+1}{2}-b}
{\frac{3}{2},\;\frac{a+1}{2}}{\,1\,},
\end{eqnarray}
and similarly for $n<0$, see \cite{mykummer}. In fact, formula
(\ref{gkummer}) is a contiguous relation between
$\hpg{2}{1}{a+n,\,b}{a-b}{-1}$, $\hpg{2}{1}{a-1,\;b}{a-b}{-1}$ and
$\hpg{2}{1}{a,\;b}{1+a-b}{-1}$, where the last two terms are
evaluated in terms of $\Gamma$-functions using Kummer's identity
(\ref{kummer}), and the coefficients are expressed as terminating
$\hpgo{3}{2}(1)$ series. % (as a bonus).

Formulas (35) and (36) in \cite{mykummer} present similar evaluations of
\[
\hpg{2}{1}{-a,\,1/2}{2a+3/2+n}{\frac{1}{4}} \quad \mbox{and} \quad
\hpg{3}{2}{a+n,\,b,\,c}{a-b,\,a-c}{1}
\]
when $n$ is an integer. These series are contiguous to the $\hpgo{2}{1}(1/4)$
and $\hpgo{3}{2}(1)$ series evaluable by Gosper's or Dixon's (respectively)
identities. Both new formulas are also three-term contiguous relations,
where two hypergeometric terms are evaluated and the coefficients are written
as terminating hypergeometric series. In both cases all series contiguous
to Gosper's $\hpgo{2}{1}(1/4)$ or well-posed $\hpgo{3}{2}(1)$ series can be
evaluated by contiguous relations, but general expressions for the
coefficients in the final three-term expression like (\ref{gkummer})
are not known.

Formulas like (\ref{gkummer}) may specialize to one-term evaluations.
%\[P(-5)=\frac{4\,(a-2)\,(a-4)\,(2\,a^2-4\,a\,b+b^2-12\,a+17\,b+12)}
%{(b-1)\,(b-2)\,(b-3)\,(b-4)}.\]
For example, $P(-5)=0$ if $2\,a^2-4\,a\,b+b^2-12\,a+17\,b+12=0$,
and under this condition formula (\ref{gkummer}) may be brought to
\begin{equation} \label{specfo2}
\hpg{2}{1}{a-5,\,b\,}{a-b}{-1}=\frac{a-b-1}{a-2\,b}\;
\frac{\Gamma(a-b-2)\;\Gamma\left(\frac{a}{2}-1\right)}
{\Gamma(a-3)\;\Gamma\left(\frac{a}{2}-b\right)}.
\end{equation}
By parameterizing the curve given by the relation between $a$ and $b$ one
gets an exotic formula, see \cite[(33)]{mykummer}. Apparently, this kind
of formula can be obtained only by using contiguous relations and
a known evaluation.
%A known way to find one-term evaluations like (\ref{specfo2}) is
%to look for hypergeometric functions which satisfy a first order
%recurrence relation, see \cite[?]{zeilb}. Such recurrence
%relations can be found by looking for two-term specializations of
%general contiguous relations. This leads to equations like
%$\bfo{Q}(k,l,m)=0$ (in the notation of the previous section).
\\

\noindent \bfo{Transformations.} General transformations of hypergeometric
series can be derived from the symmetries of their contiguous relations. For
example, all terms in the relations between the 24 Kummer's $\hpgo{2}{1}$
functions (see \cite[2.9]{bateman}) satisfy not only the same hypergeometric
differential equation, but also the same contiguous relations with the same
shifts in the parameters $a$, $b$ and $c$. To show this statement, one can
check that the two functions
\begin{equation}
u_1=\hpg{2}{1}{a,\,b\,}{c}{z},\quad u_2=\frac{\Gamma(c)\,\Gamma(c-a-b)}
{\Gamma(c-a)\,\Gamma(c-b)}\,\hpg{2}{1}{a,\,b}{a\!+\!b\!+\!1\!-\!c}{1\!-\!z}
\end{equation}
satisfy basic contiguous relations (\ref{contigaa})--(\ref{contigab}). Then
the functions
\begin{eqnarray}
u_3&=&\frac{\Gamma(c)\,\Gamma(b-a)}{\Gamma(c-a)\,\Gamma(b)}\;(-z)^{-a}\;
\hpg{2}{1}{a,\,a+1-c\,}{a+1-b}{\frac{1}{z}},\\
u_4&=&\frac{\Gamma(c)\,\Gamma(a-b)}{\Gamma(c-b)\,\Gamma(a)}\;(-z)^{-b}\;
\hpg{2}{1}{b+1-c,\,b\,}{b+1-a}{\frac{1}{z}},\\
u_5&=&\frac{\Gamma(c)\,\Gamma(1\!-\!a)\,\Gamma(1\!-\!b)}
{\Gamma(2\!-\!c)\,\Gamma(c\!-\!b)\,\Gamma(c\!-\!a)}\,z^{1-c}\,
\hpg{2}{1}{a\!+\!1\!-\!c,\,b\!+\!1\!-\!c\,}{2-c}{z},\\
u_6&=&\frac{\Gamma(c)\,\Gamma(a+b-c)}{\Gamma(a)\,\Gamma(b)}\,(1\!-\!z)^{c-a-b}
\,\hpg{2}{1}{c-a,\,c-b\,}{c\!+\!1\!-\!a\!-\!b}{1\!-\!z},
\end{eqnarray}
satisfy the same contiguous relations as well, since the
coefficients in the relations
\begin{eqnarray}
\label{kum24a} \!\!u_6=u_1\!-u_2, & & \hspace{-14pt}\frac{\bfo{e}^{-i\pi
a}\sin(\pi(b\!-\!a))}{\sin(\pi b)}\,u_3=u_1\!-
\frac{\bfo{e}^{i\pi(c-a)}\sin(\pi(a\!+\!b\!-\!c))}{\sin(\pi b)}u_6,\\
\label{kum24b} \!\!u_4=u_1\!-u_3, &  & \hspace{69pt}u_5=u_1-\frac{\sin(\pi
c)\sin(\pi(a\!+\!b\!-\!c))}{\sin(\pi a)\,\sin(\pi b)}\,u_6,
\end{eqnarray}
%u_1=u_5+\frac{\sin(\pi a)\,\sin(\pi b)}{\sin(\pi c)\,\sin(\pi(a+b-c))}\,u_6,
%u_1=\frac{\bfo{e}^{i\pi a}\,\sin(\pi b)}{\sin(\pi(b-a))}\,u_3+
%\frac{\bfo{e}^{i\pi(c-a)}\,\sin(\pi b)}{\sin(\pi(a+b-c))}\,u_6.
are constants with respect to integral shifts in the parameters $a$, $b$,
$c$. Other 18 hypergeometric functions are alternative representations of
$u_1,\ldots,u_6$, see \cite[2.9]{bateman}. All relations between the 24
Kummer's functions are generated by (\ref{kum24a})--(\ref{kum24b}), so the
statement follows. Mind that different expressions of $u_i$'s as
hypergeometric functions may have different arguments.

In principle, relations (\ref{kum24a})--(\ref{kum24b}) can be found by
showing that the three terms satisfy the same second order recurrence
relation (with respect to integral shifts in some parameter), and comparing
their asymptotics as the corresponding parameter approaches $\infty$ and
$-\infty$. However, this approach may be cumbersome.  Quadratic or higher
order algebraic transformations can also be found in this way.

An interesting question is whether the symmetries of contiguous relations
can classify all identities between hypergeometric series. For $\hpgo{2}{1}$
series this boils down to the symmetries of ${\bf Q}(k,l,m)$. In the
remainder of this section, we demonstrate a way to investigate the
symmetries of recurrence relations (a special case of contiguous relations)
of hypergeometric series.
\\

\noindent \bfo{Symmetries of recurrence relations.} To obtain recurrence
relations for some hypergeometric series one can introduce the discrete
parameter $n$ by replacing one or several continuous parameters as $a\mapsto
a+2n$ (and similarly) in different ways. Then a recurrence is a contiguous
relation between sufficiently many hypergeometric functions with successive
$n$. We are interested in situations when two hypergeometric series (with a
common discrete parameter $n$) satisfy the same recurrence relation, perhaps
after multiplying the functions by a hypergeometric term (that is, a
solution of a first order recurrence relation with coefficients rational in
$n$, see \cite[5.1]{zeilb}). Therefore we call two hypergeometric functions
with discrete parameter $n$ {\em equivalent} if they differ by such a
factor.

More specifically, suppose that a hypergeometric series satisfies
a second order recurrence relation (for $n\ge 0$)
\begin{equation}
A(n)\,S(n+1)+B(n)\,S(n)+C(n)\,S(n-1)=0,
\end{equation}
with $A(n)$, $B(n)$ and $C(n)$ being rational functions in $n$.
For convenience, we assume that $S(n)$ does not satisfy a first
order recurrence relation, and that $B(n)\neq 0$ for all $n\ge 0$.
Then $S(n)$ is equivalent to:
\[
Z(n)=(-1)^n\frac{A(0)\ldots A(n-1)}{B(0)\ldots B(n-1)} S(n).
\]
The sequence $Z(n)$ satisfies the recurrence:
\begin{equation} \label{seqzrec}
Z(n+1)-Z(n)+\frac{C(n)\,A(n-1)}{B(n)\,B(n-1)}\,Z(n-1)=0.
\end{equation}
In fact, $Z(n)$ is the unique sequence equivalent to $S(n)$ and
satisfying a recurrence relation of form $Z(n+1)-Z(n)+H(n)
Z(n-1)=0$. (Note that the second order recurrence for $Z(n)$ is
unique up to a factor, and that we want to keep the coefficients
of two terms.) It follows that an equivalence class of
hypergeometric functions is determined by the rational function
\begin{equation} \label{eqclassf}
\frac{B(n)\,B(n-1)}{C(n)\,A(n-1)}.
\end{equation}

Hypergeometric series with a discrete parameter $n$ can be
classified by their equivalence class function. One can compute
this function for various types of hypergeometric functions and
different ways of introducing the discrete parameter $n$. Then one
can look for the cases when different equivalence class functions
are equal (perhaps under some relations between their continuous
parameters).

For example, the equivalence class function for
$\hpg{0}{1}{-}{\hat{c}+n}{\hat{z}}$ is the polynomial
$-(n\!+\!\hat{c}\!-\!1)(n\!+\!\hat{c}\!-\!2)/\hat{z}$. This is
also the equivalence class function for
$\hpg{0}{1}{-}{2-\hat{c}-n}{\hat{z}}$. The equivalence class
function for $\hpg{1}{1}{a+n}{c+2n}{z}$ is
\[
-\frac{(2n\!+\!c\!-\!1)(2n\!+\!c\!-\!3)(4n^2\!+\!4cn\!-\!4n\!+\!c^2\!+
\!2az\!-\!cz\!-\!2c)(4n^2\!+\!4cn\!-\!12n\!+\ldots)}
{z^2\,(n+a-1)\,(n+c-a-1)\,(2n+c)\,(2n+c-4)}.
\]
When $c=2a$, this rational function is also a polynomial of degree two,
namely $-4(2n\!+\!2a\!-\!1)(2n\!+\!2a\!-\!3)/z^2$. By comparing these two
polynomials and corresponding transformations of the recurrences to
normalized form (\ref{seqzrec}), we conclude that there must be a linear
relation between the three functions %$$,
\begin{eqnarray*}
\hpg{1}{1}{a+n}{\!2a+2n\!}{\!4z},\qquad\hpg{0}{1}{-}{\!a\!+\!\frac{1}{2}\!+\!n}{z^2}
\hspace{3pt}\\ %\quad
\mbox{and}\qquad
\frac{\Gamma(a\!+\!\frac{1}{2}\!+\!n)}{\Gamma(\frac{3}{2}\!-\!a\!-\!n)}
\;z^{1-2a-2n}\;\hpg{0}{1}{-}{\frac{3}{2}\!-\!a\!-\!n}{z^2},
\end{eqnarray*}
where coefficients are constants with respect to the shift $n\to n\!+\!1$.
Note that the first two functions are bounded as $n\to\infty$, whereas the
last one is unbounded (for fixed generic $a$ and $z$). Hence the coefficient
of the last function is zero. By comparing the limits of the first two
functions one concludes that
$\hpg{1}{1}{a+n}{2a+2n}{4z}=\bfo{e}^{2z}\,\hpg{0}{1}{-}{a+n+\frac{1}{2}}{z^2}$;
see also \cite[(4.1.12)]{specfaar}.
%For each rational function one can compute the resultant of its
%numerator and denominator with respect to $n$ in order to
%investigate cases when the rational function can be simplified to
%a quotient of polynomials of smaller degree in $n$. Then different
%equivalence class functions with the same degrees of the numerator
%and the denominator can be compared.

The author computed a number of equivalence class functions for
$\hpgo{r}{1}$ (with $r=0,1,2$) and $\hpgo{3}{2}(1)$ functions. They do not
imply new identities, except plenty of exotic two-term identities such as
(\ref{specfo2}), and some interesting consequences of known identities.
Straightforward computations are too cumbersome when the degree of
numerators or denominators of equivalence class functions exceeds 10.
However, this method generalizes easily to $q$-hypergeometric functions. A
more intelligent consideration of symmetries of recurrence (and contiguous)
relations could be helpful in finding relations between hypergeometric
functions when one cannot use symmetries of differential equations for them,
say when the variable $z$ is specialized, or when the parameters and $z$ are
related.
%Recurrence relations help to prove identities where the variable
%$z$ is specialized, or the parameters and $z$ are related.

%The terms in contiguous relations usually do not satisfy the same recurrence
%relations.

\section{Proof of Theorem \ref{spcontig}}
\label{thproofs}

To show the existence of (\ref{contklm}), observe that contiguous relations
(\ref{contigaa})--(\ref{contigab}) express $\ghpg{a-1,\,b\!}{c}$,
$\ghpg{a,\,b}{c-1}$, $\ghpg{a,\,b+1\!}{c}$ in terms of $\ghpg{a,\,b}{c}$ and
$\ghpg{a+1,\,b}{c}$. There are similar expressions for $\ghpg{a,\,b}{c+1}$
and $\ghpg{a,\,b-1\!}{c}$. Using shifted versions of these relations one can
express $\ghpg{\!a+k,\,b+l\!}{c+m\!}$ linearly in hypergeometric functions
without shifts in the parameters $b$ and $c$, and then leave only terms
$\ghpg{a,\,b}{c}$ and $\ghpg{a+1,\,b}{c}$.

If expression (\ref{contklm}) is not unique, then
$\ghpg{a+1,\,b}{c}\Big/\ghpg{a,\,b}{c}$ is a rational function in the
parameters, just as $\ghpg{a+1,\,b+1}{c}\Big/\ghpg{a+1,\,b}{c}$ (by the
symmetry of the upper parameters). Then
$\hpg{2}{1}{a+1,\,b+1}{1+a-b}{\!-1}\Big/\hpg{2}{1}{a,\,b}{1+a-b}{\!-1}$ must
be a rational function in $a$ and $b$. But this function has an unbounded
set of poles according to Kummer's identity (\ref{kummer}). Hence a
contradiction.

The uniqueness of (\ref{contklm}) implies that $\bfo{P}(k,l,m)$ and
$\bfo{Q}(k,l,m)$ satisfy the same contiguous relations as $\hpgo{2}{1}$
series; check the triple substitution of (\ref{contklm}) into
(\ref{gencontig}).

The ``initial" conditions (\ref{initcond}) are obvious.

To prove (\ref{pinq}), apply (\ref{contklm}) to both sides of
$\ghpg{a+k,\,b}{c}=\shif{a}\,\ghpg{a+k-1,\,b}{c}$, and use
(\ref{contigaa}) to eliminate $\ghpg{a-1,\,b}{c}$ on the
right-hand side.

Formula (\ref{qargplus}) is obtained after several applications of
(\ref{contklm}) to $\ghpg{a+k+k',\,b+l+l'}{c+m+m'}$. The first
intermediate step is:
\begin{eqnarray*}
\ghpg{a+k+k',\,b+l+l'}{c+m+m'}&\!\!=\!\!&\left(\shif{a}^k\shif{b}^l\shif{c}^m\,
\bfo{P}(k',l',m')\right)\,\ghpg{a\!+\!k,b\!+\!l}{c+m}\\
&&\hspace{-10pt}+\left(\shif{a}^k\shif{b}^l\shif{c}^m\,
\bfo{Q}(k',l',m')\right)\,\ghpg{a\!+\!k\!+\!1,b\!+\!l}{c+m}.\hspace{-4pt}
\end{eqnarray*}
Now apply (\ref{contklm}) to get terms with $\ghpg{a,\,b}{c}$ and
$\ghpg{a+1,\,b}{c}$ only, compare the terms to
$\ghpg{a+1,\,b}{c}$, and eventually use (\ref{pinq}) once to
replace $\bfo{P}(k',l',m')$.

For a proof of the last two formulas, let us introduce:
\begin{equation}
W_{p,q,r}(k,l,m):=\left|\begin{array}{cc}
\bfo{P}(k,\,l,\,m) & \bfo{P}(k+p,\,l+q,\,m+r) \\
\bfo{Q}(k,\,l,\,m) & \bfo{Q}(k+p,\,l+q,\,m+r) \end{array}\right|.
\end{equation}
We assume that $p$, $q$, $r$ (just as $k$, $l$, $m$) are integers.
\begin{lemma} \label{wpqrlemma}
The following properties of the $W$-symbol hold.
\begin{enumerate}
\item $W_{0,0,0}(k,l,m)=0$. \item $W_{p,q,r}(0,0,0)=\bfo{Q}(p,q,r)$. \item
For fixed $k,l,m$ the determinants $W_{p,q,r}$ satisfy the contiguous
relations of $\hpgo{2}{1}$ functions. More precisely, if {\rm
(\ref{gencontig})} is a contiguous relation, then $\sum_{n=1}^3 \left(
\shif{a}^{k+p}\,\shif{b}^{l+q}\,\shif{c}^{m+r}\,f_n\right)\;
W_{p+\alpha_n,q+\beta_n,r+\gamma_n}(k,\,l,\,m)=0$. \item $W_{1,0,0}$
satisfies the first order recurrence relations
\begin{eqnarray*}
W_{1,0,0}(k+1,l,m)&=&\left( \shif{a}^{k+1}\,\shif{c}^m\,
\,\frac{a-c}{a\,(1-z)} \right)\, W_{1,0,0}(k,l,m),\\
W_{1,0,0}(k,l+1,m)&=&\left( \shif{b}^l\;\shif{c}^m\,
\,\frac{b-c+1}{b\,(1-z)}\right)\,W_{1,0,0}(k,l,m),\\
W_{1,0,0}(k,l,m+1)&=&\left( \shif{a}^k\,\shif{b}^l\,\shif{c}^m
\,\frac{c^2\,(z-1)}{(c-a)\,(c-b)\,z}\right)\,W_{1,0,0}(k,l,m).
\end{eqnarray*}
\end{enumerate}
\end{lemma}
\proof The first two properties are straightforward.

If (\ref{gencontig}) is a contiguous relation, then
\begin{eqnarray*}
\sum_{n=1}^3 \left(
\shif{a}^{k+p}\,\shif{b}^{l+q}\,\shif{c}^{m+r}\,f_n \right)\;
W_{p+\alpha_n,q+\beta_n,r+\gamma_n}(k,\,l,\,m) \ = \\
\left|\begin{array}{cc} \bfo{P}(k,\,l,\,m) & \sum_{n=1}^3 \left(
\shif{a}^{k+p}\shif{b}^{l+q}\shif{c}^{m+r}\,f_n\right)\,\bfo{P}
(k\!+\!p\!+\!\alpha_n,l\!+\!q\!+\!\beta_n, m\!+\!r\!+\!\gamma_n)
\\ \bfo{Q}(k,\,l,\,m) & \sum_{n=1}^3 \left(
\shif{a}^{k+p}\shif{b}^{l+q}\shif{c}^{m+r}\,f_n\right)\,\bfo{Q}
(k\!+\!p\!+\!\alpha_n,l\!+\!q\!+\!\beta_n,
m\!+\!r\!+\!\gamma_n) \end{array}\right| \\
=\ \left|\begin{array}{cc} \bfo{P}(k,l,m) & 0 \\
\bfo{Q}(k,l,m) & 0 \end{array}\right| \ =\ 0.
\end{eqnarray*}

The recurrence relation for $W_{1,0,0}$ with respect to $k$
follows from contiguous relation (\ref{contigaa}):
\begin{eqnarray*}
W_{1,0,0}(k\!+\!1,l,m) & = & \left|\begin{array}{cc}
\bfo{P}(k+1,l,m)
& \emphi{A}\,\bfo{P}(k,l,m)+\emphi{B}\,\bfo{P}(k\!+\!1,l,m) \\
\bfo{Q}(k+1,l,m) & \emphi{A}\,\bfo{Q}(k,l,m)
+\emphi{B}\,\bfo{Q}(k\!+\!1,l,m) \end{array}\right| \\
& = & -\emphi{A}\,W_{1,0,0}(k,l,m), \quad\mbox{with }
\emphi{A}=\shif{a}^{k+1}\shif{c}^m\,\frac{c-a}{a\,(1-z)}.
%\frac{c-a+m-k-1}{(a+k+1)\,(1-z)}.
\end{eqnarray*}

To prove the recurrence relation for $W_{1,0,0}$ with respect to
$l$ we use contiguous relations (\ref{contigab}) and
\cite[2.8.(36)]{bateman} to derive
\begin{eqnarray*}
\bfo{P}(k,l+1,m) & = & A\,\bfo{P}(k,l,m)+B\,\bfo{P}(k+1,l,m), \\
\bfo{P}(k+1,l+1,m) & = & C\,\bfo{P}(k,l,m)+D\,\bfo{P}(k+1,l,m),
\end{eqnarray*}
with $A=\shif{a}^k\shif{b}^l\,\frac{b-a}{b}$,
$B=\shif{a}^k\shif{b}^l\,\frac{a}{b}$,
$C=\shif{a}^k\shif{b}^l\shif{c}^m\,\frac{c-a-1}{b\,(1-z)}$,
$D=\shif{a}^k\shif{b}^l\shif{c}^m\,\frac{a+b-c+1}{b\,(1-z)}$.
There are similar expressions for $\bfo{Q}(k,l+1,m)$ and
$\bfo{Q}(k+1,l+1,m)$. Then
\begin{eqnarray*}
W_{1,0,0}(k,\,l+1,\,m) \ = \\
\left|\begin{array}{cc} A\,\bfo{P}(k,l,m)+B\,\bfo{P}(k\!+\!1,l,m)
& C\,\bfo{P}(k,l,m)+D\,\bfo{P}(k\!+\!1,l,m) \\
A\,\bfo{Q}(k,l,m)+B\,\bfo{Q}(k\!+\!1,l,m) & C\,\bfo{Q}(k,l,m)
+D\,\bfo{Q}(k\!+\!1,l,m) \end{array}\right| \ = \\
(A\,D-B\,C)\;W_{1,0,0}(k,\,l,\,m),
\end{eqnarray*}
which is the required relation. The recurrence with respect to $m$ can be
proved similarly, using formulas (34) and (38) in \cite[2.8]{bateman}. This
completes the proof of the lemma.
%$\frac{c\,(cz-bz-a)}{(c-a)\,(c-b)\,z}$, $\frac{a\,c\,(1-z)}{(c-a)\,(c-b)\,z}$,
%$\frac{c}{(c-b)\,z}$ and $\frac{c\,(z-1)}{(c-b)\,z}$.
\QED

The recurrence relations for $W_{1,0,0}(k,l,m)$ imply that
\begin{equation} \label{woneklm}
W_{1,0,0}(k,l,m)=\frac{(c)_m\;(c)_m\;z^{-m}\,(z-1)^{m-k-l}}
{(a\!+\!1)_k\,(b)_l\,(c\!-\!a)_{m-k}\,(c\!-\!b)_{m-l}}.
\end{equation}
The recurrence relations for $W_{p,q,r}$ and an initial condition
{\em (i)} imply that
\begin{equation} \label{wpqrklm}
W_{p,q,r}(k,l,m)=W_{1,0,0}(k,l,m)\cdot\shif{a}^k\,\shif{b}^l\,\shif{c}^m\,
\bfo{Q}(p,q,r).
\end{equation}
Formula (\ref{pqdet}) follows from these two equations.

To prove (\ref{minklm}) note that the determinant in (\ref{pqdet})
can be expressed in two ways as the symbol $W$. This gives:
\begin{equation}
W_{k-k',l-l',m-m'}(k',l',m')=-W_{k'-k,l'-l,m'-m}(k,l,m).
\end{equation}
By setting $k'=l'=m'=0$ and changing the sign of $k$, $l$ and $m$ one obtains:
\begin{equation}
W_{-k,-l,-m}(0,0,0)=-W_{k,l,m}(-k,-l,-m).
\end{equation}
Statement {\em (ii)} of lemma \ref{wpqrlemma}, formulas
(\ref{wpqrklm}) and (\ref{woneklm}), and transformation of
Pochhammer symbols give the last formula of Theorem
\ref{spcontig}. \QED

\bibliographystyle{alpha}
\bibliography{../hypergeometric}

\begin{thebibliography}{PWZ96}

\bibitem[AAR99]{specfaar}
G.E. Andrews, R.~Askey, and R.~Roy.
\newblock {\em Special Functions}.
\newblock Cambridge Univ. Press, Cambridge, 1999.

\bibitem[Erd53]{bateman}
A.~Erd\'elyi, editor.
\newblock {\em Higher Transcendental Functions}, volume~I.
\newblock McGraw-Hill Book Company, New-York, 1953.

\bibitem[Knu71]{knuthii}
D.E. Knuth.
\newblock {\em The art of computer programming}, volume II, Seminumerical
  Algorithms.
\newblock Addison-Wesley, Reading, Mass., 1971.

\bibitem[Pau01]{paulecont}
P.~Paule.
\newblock Contiguous relations and creative telescoping.
\newblock Technical report, RISC, Austria, 2001.

\bibitem[PWZ96]{zeilb}
M.~Petkov\v{s}ek, H.S. Wilf, and D.~Zeilberger.
\newblock {\em A=B}.
\newblock A. K. Peters, Wellesley, Massachusetts, 1996.

\bibitem[Rai60]{rainvill}
E.D. Rainville.
\newblock {\em Special Functions}.
\newblock The MacMillan Company, New York, 1960.

\bibitem[Tak89]{takayamcontig}
N.~Takayama.
\newblock Gr\"{o}bner basis and the problem of contiguous relations.
\newblock {\em Japan Journal of Applied Mathematics}, 6:147--160, 1989.

\bibitem[Vid03]{mykummer}
R.~Vid\=unas.
\newblock A generalization of {K}ummer's identity.
\newblock {\em Rocky Mount. J.}, 32, 2003.
\newblock also available at {\sf http://arXiv.org/abs/math.CA/0005095}.

\bibitem[VK00]{koornvid}
R.~Vid\=unas and T.H. Koornwinder.
\newblock Web-page of the {NWO} project "{A}lgorithmic methods for special
  functions by computer algebra".
\newblock {\sf http://www.science.uva.nl/$\sim$thk/specfun/compalg.html}, 2000.

\end{thebibliography}

\end{document}